\documentclass[12pt]{amsart}
\usepackage{bigints}

\usepackage{graphicx}

\usepackage{pdfsync}
\usepackage{dsfont}
\usepackage{mathrsfs}

\usepackage{mathptmx}
\usepackage{helvet}
\usepackage{courier}
\usepackage{multicol}
\usepackage{footmisc}
\usepackage{enumerate}

\usepackage{float}
\restylefloat{table}

\usepackage{pgf,tikz}
\usepackage{pgfplots}
\usetikzlibrary{automata}
\usetikzlibrary{arrows}
\usetikzlibrary{positioning}

\newtheorem{Th}{Theorem}[section]

\newtheorem{Ex}[Th]{Example}

\newtheorem{Rem}[Th]{Remark}

\usepackage{bm}

\usepackage{tensor}
\usepackage{comment}

\pgfplotsset{compat=1.10}

\newcommand{\Exp}{\mathds{E}}

\newcommand{\Prob}{\mathds{P}}

\newcommand{\vect}[1]{\vec{#1}}
\renewcommand{\vect}[1]{\boldsymbol{#1}}

\newcommand{\mat}[1]{\boldsymbol{#1}}

\renewcommand{\dfrac}[2]{\frac{\displaystyle {#1}}{\displaystyle {#2}}}

\renewcommand{\matrix}[2][ccccccccccccccccccccc]{\left(\begin{array}{#1}#2
      \\ \end{array} \right)}

\DeclareMathOperator*{\argmax}{arg\,max}

\usepackage[top=1.75cm,bottom=1.75cm,left=2.5cm,right=2.5cm]{geometry}

\title{Fitting phase--type scale mixtures to heavy--tailed data and distributions}

 \author[M. Bladt]{Mogens Bladt}
 \address{Institute for Applied Mathematics and Systems \\
 National University of Mexico \\
 A.P. 20-726 \\
 01000 Mexico, D.F.\\
 Mexico}
 \email{bladt@sigma.iimas.unam.mx}

\author[L. Rojas]{Leonardo Rojas-Nandayapa}
\address{School of mathematics and physics\\
The University of Queensland \\
St. Lucia 4072, Brisbane \\
Australia}
\email{l.rojas@uq.edu.au}
\begin{document}

\allowdisplaybreaks
\begin{abstract}

{We consider the fitting of heavy tailed data and distribution with a special attention to distributions with a non--standard shape in the ``body'' of the distribution. To this end we consider a dense class of heavy tailed distributions introduced in \cite{bladt-nielsen-samorodnitsky:2015}, employing an EM algorithm for the the maximum likelihood estimates of its parameters. 
We present methods for fitting to observed data, histograms, 
censored data, as well as to theoretical distributions. 
Numerical examples are provided with simulated data and a benchmark reinsurance dataset.  
We empirically demonstrate that our model can provide excellent fits to heavy--tailed data/distributions
with minimal assumptions. 
}

\smallskip
\noindent \textbf{Keywords.} Statistical inference, heavy--tailed, phase--type, scale mixtures, approximating distributions,
  EM algorithm.
\end{abstract}

\maketitle

\section{Introduction}
In this paper we consider the maximum likelihood estimation for a dense class of nonnegative heavy--tailed tailed distributions, referred to as 
scale mixtures of phase--type distribuitions (NPH), which was defined in \cite{bladt-nielsen-samorodnitsky:2015}. 
{Distributions in the NPH class
allow for the simultaneous modelling of the ``body'' and the ``tail'' of general distributions which are assumed to be 
absolutely continuous and nonnegative while their ``tails''  are assumed to belong to some general class of heavy tailed distributions, 
like for instance Regularly Varying or Weibullian. These very general assumptions allows us to fit
heavy--tailed distributions which may look distinctively different from distributions usually found in catalogues.}

 Apart from providing an adequate description of  data,  distributions from NPH can also be seen as infinite--dimensional phase--type distributions, which to all intents and purposes are as tractable
as their finite--dimensional counterparts. 
Much of the machinery available for finite--dimensional phase--type distributions is also applicable to the extended class.
Algorithms are for example available for the exact calculations of properties related to renewal theory, random walks (ladder processes) and ruin probabilities
(see \cite{bladt-nielsen-samorodnitsky:2015} for details). 

The maximum likelihood estimation will be carried out employing an EM algorithm 
similarly as for finite--dimensional phase--type distributions \cite{ANO:1996}. The main challenge we face is the algorithmic implementation 
resulting from the extension to infinite dimensions since we cannot make a pre--fixed cut--off in the number of dimensions as this would be equivalent to a light--tailed estimation. We shall see, however, that it is possible to reduce the formulas to simple (one--dimensional) infinite series which can be numerically evaluated to a specified degree of precision. 
We are
also able to simplify certain expressions of the original paper (\cite{ANO:1996}), obtaining explicit expressions which allow for an improved numerical performance.

The rest of the paper is organized as follows. In Section \ref{sec:estimation} we develop the main algorithm for estimating independent and identically distributed data sampled from an NPH distribution. We show that the algorithm essentially uses the empirical cumulative distribution function, and that a significant increase in speed may be obtained by grouping the data into intervals on parts of the support and considering the corresponding histograms. In Section \ref{sec:censoring} we adjust the algorithm to cope with the presence of left-, right- or interval censored data. Section \ref{sec:KL} provides an EM--algorithm for adjusting an NPH distribution to a theoretical distribution function $F$, which in turn is equivalent to 
finding the distribution in the NPH class with a given order which minimizes the Kullback--Leibler distance to $F$.  
In section \ref{sec:numerical} we provide some numerical examples from both simulated data, real data and a fits to theoretical distribution.
In there, we highlight some minimal assumptions made for efficiently fitting general Regularly
Varying and Weibullian distributions.

\section{Phase--type distributions and the NPH class}

Before proceeding with a more detailed account, we first provide some background on phase--type distributions and the extended class $\mbox{NPH}$ of scale
mixtures of phase--type distributions.

The class $\mbox{PH}$  of phase--type distributions consists of distributions of
(random) times until a finite state Markov jump process exits a set of transient states. This can be made precise by letting $E=\{1,2,...,p,p+1\}$ denote the state--space of a Markov jump process $\{ X_t\}_{t\geq 0}$, where states $1,2,...,p$ are transient and $p+1$ is absorbing. The intensity matrix for $\{ X_t\}_{t\geq 0}$ can then be written on the form
\[  \mat{\Lambda}  = \matrix{\mat{T} & \vect{t} \\ \vect{0} & 0} , \]
where $\mat{T}$ is a $p\times p$ sub--intensity matrix, $\vect{t}$ is a $p$--dimensional column vector and $\vect{0}$ the $p$--dimensional row vector of zeroes. We follow the convention that matrices are written in capital bold and their elements with the corresponding minuscule  (like $\mat{A}=\{ a_{ij}\}$), bold minuscule Greek letters (like $\vect{\alpha}$) are row vectors while bold minuscule Roman letters (like $\vect{t})$ are column vectors. Their dimensions are usually clear from the context and left unspecified unless necessary. 
Let $\alpha_i=\Prob (X_0=i)$, $\sum_{i=1}^p \alpha_i=1$ and $\vect{\alpha}=(\alpha_1,...,\alpha_p)$. Then we say that 
\[ \tau = \inf\{ s>0 : X_s = p+1 \} \]
has a phase--type distribution with representation $\mbox{PH}_p(\vect{\alpha},\mat{T})$. 
Since $\vect{t}=-\mat{T}\vect{e}$, where $\vect{e}$ denotes the column vector of ones, the distribution of $\tau$ is fully specified in terms of $\vect{\alpha}$ and $\mat{T}$. For further background on phase--type distributions we refer e.g. to \cite{memap}, \cite{scholarometermbladt6}, \cite{Latouche:1987uw} or \cite{Neuts:1981vg}.

The class of phase--type distributions is widely used in the area of Applied Probability, where they may often provide exact (or even explicit) solutions in complex
stochastic models. This is for example the case for ruin probabilities in risk theory or waiting time distributions for queues. Any distribution with support on the positive real numbers may be approximated arbitrarily close by a phase--type distribution. In spite of this denseness property,  phase--type distributions are all light tailed and consequently any approximating (finite--dimensional) phase--type distribution will not be able to capture a possible heavy tailed behaviour. 

In \cite{bladt-nielsen-samorodnitsky:2015} the dense class, $\mbox{NPH}$, of genuinely heavy tailed distributions is proposed in terms of infinite--dimensional phase--type distributions with a finite number of parameters. The idea is very simple and goes a follows. Let $N$ be a discrete random variable with support on $\{ s_i: {i\in\mathds{N}},s_i>0\}$ and distribution $\vect{\pi}=\{ \pi_i \}_{i\geq 1}$ where $\pi_i = \Prob (N=s_i)$. If $N$, for example, is a discretization of a continuous distribution $G$ at step length $\Delta$, then  $s_i= i \Delta$ and $\pi_i =G(i\Delta)-G((i-1)\Delta)$, $i=1,2,\dotsc $. Let $\tau \sim \mbox{PH}_p(\vect{\alpha},\mat{T})$ be independent of $N$. When sampling from $Y$, we first draw an index $i$ (referred to as the {\it level}) from $\vect{\pi}$ and then a phase--type random variable from $\mbox{PH}_p(\vect{\alpha},\mat{T}/s_i)$.
Hence, the random variable $Y=N\tau$ may be seen as a scale mixture of phase--type distributions,
so we often refer to $N$ as the scaling random variable and its distribution $\vect{\pi}$ as the scaling distribution of the NPH.
The density of $Y$ can written as
\[  f_Y(y) = \sum_{i=1}^\infty \pi_i \vect{\alpha} e^{\mat{T} y/s_i}\vect{t}/s_i ,\qquad y>0. \]
The distribution of $Y$ can also be seen as an infinite dimensional phase--type distribution since we can write
\[   f_Y(y)= \left( \vect{\pi}\otimes \vect{\alpha}\right) e^{\mat{\Gamma} y } \vect{\gamma} ,  \]
where
\[  \mat{\Gamma} =\matrix{\mat{T}_1 & \mat{0} & \mat{0} & \mat{0} &  \cdots \\
\mat{0} & \mat{T}_2 & \mat{0} & \mat{0} & \cdots \\
\mat{0} & \mat{0} & \mat{T}_3 & \mat{0} & \cdots  \\
\mat{0} & \mat{0} & \mat{0} & \mat{T}_4 & \cdots \\
\vdots & \vdots & \vdots & \vdots & \ddots }= 
\matrix{\mat{T}/s_1 & \mat{0} & \mat{0} & \mat{0} &  \cdots \\
\mat{0} & \mat{T}/s_2 & \mat{0} & \mat{0} & \cdots \\
\mat{0} & \mat{0} & \mat{T}/s_3 & \mat{0} & \cdots  \\
\mat{0} & \mat{0} & \mat{0} & \mat{T}/s_4 & \cdots \\
\vdots & \vdots & \vdots & \vdots & \ddots }  \] 
and $\vect{\gamma} = -\mat{\Gamma}\vect{e}$. Here $\vect{e}$ is now the infinite--dimensional column vector of ones and we notice that the exponential of the infinite dimensional matrix  $\mat{\Gamma}$ is well defined since it is a {\it bounded operator} (as the sequence $\{ s_i:s_i>0\}$ is bounded away from zero). We also let $\vect{t}_i = - \mat{T}_i\vect{e}$. We shall write
\[ Y \sim \mbox{NPH}_p (\vect{\pi}, \vect{\alpha},\mat{T}) .  \]
While the support of $\vect{\pi}$ is important, we will not denote it explicitly in the parametrisation 
of $Y$.  

A key feature of the NPH class is that it contains a rich variety of  heavy--tailed distributions. 
For instance, if $\vect{\pi}$ has unbounded support, then $Y$ has
a heavy--tailed distribution \cite[cf.][]{RojasXie}. In the regularly varying case, 
the tail of $Y$ greatly follows that of $N$.  Breiman's lemma implies that if
$N$ has a regularly varying with tail index $-\alpha<0$, then the tail of $Y$ will also be regularly varying with the
same index.  More precisely,
\begin{equation*}
 \Prob(Y>t)=\Exp[\tau^\alpha]\Prob(N>t)(1+o(1)),\qquad t\to\infty.
\end{equation*} 
A similar feature occurs for the
the class of Weibullian distributions  \cite{ArendarczykDebicki2011}, which is 
defined as the collection of nonnegative distributions having survival function
\begin{equation}
 \overline F(x)= x^\delta e^{-(\lambda x)^p}\,(C+\omicron(1)),\qquad x>0, \lambda,p>0,\delta\in \mathbb{R}.
\end{equation}
If the parameter $p\in(0,1)$, then the distribution is heavy--tailed (and light--tailed otherwise). 
Since a PH distribution is Weibullian with parameter $p=1$, then Lemma 2.1 of  \cite{ArendarczykDebicki2011} implies 
that if the scaling distribution $\vect\pi$ is Weibullian with parameter $p_1$ then
$Y\sim\mbox{NPH}(\vect\pi,\vect\alpha,\mat T)$ is Weibullian with parameter $p_1/(1+p_1)$.

\section{Estimation}\label{sec:estimation}


{In this section we address the problem of fitting an NPH distribution to a data set.
We assume that $y_1,y_2,...,y_M$ forms an i.i.d.\ data set sampled from $\mbox{NPH}_p(\pi (\vect{\theta}),\vect{\alpha},\mat{T})$, 
where $\vect\pi(\vect{\theta})$ is some parametric distribution describing the distribution of $N$. We assume that the support for $\vect{\pi}(\vect{\theta})$ does not depend on $\vect{\theta}$.
We shall estimate the parameters $\vect{\theta},\vect{\alpha}$ and $\mat{T}$. }

Hence, with probability $\vect{\pi}(\vect{\theta})_i$, $y_n$ is the realization of the $i$--th level phase--type distribution 
$\mbox{PH}_p(\vect{\alpha},\mat{T}/s_i)$, but both the
level $i$ and the actual realization of the underlying Markov jump process
are not observable.  Thus, we resort to the EM algorithm for approximating the maximum likelihood estimators. 
To this end we first attend the corresponding estimation problem for complete data. 

Assume that apart from $y_1,...,y_M$ we have also observed all levels and all underlying Markov jump process. 
Let $I_n$
denote the level of the phase--type distribution of the $n$'th path,
and let
\[ L^i = \sum_{n=1}^M 1\{ I_n=i\}  \]
denote the number of $i$--level processes in the data. 
Next consider the Markov jump
process $\{ J_u^{(n)}\}_{u\geq 0}$ underlying the $n$'th phase--type distribution (which generates the data $y_n$) 
and let
\[ B_k^i=\sum_{n=1}^M 1\{ J_0^{(n)}=k, I_n=i\}  \]
be the number of $i$--level processes that are initiated in state $k$.
Define
\[ Z_k^i =\sum_{n=1}^M \int_0^{y_n} 1\{ J_u^{(n)}=k, I_n=i\}du ,  \]
which is the total time all underlying $i$--level Markov jump processes spend 
in state $k$ and let $N_{k\ell}^i$ denote the total number of jumps from state
$k$ to $\ell$ within all $i$--level Markov jump processes. Finally, let $N_k^i$
be the number of $i$--level processes that exit to the absorbing state 
from state $k$. 

Then the complete data
likelihood is easily seen to be (see e.g. \cite{ANO:1996} or \cite{memap} for further comments on this)
\[ L_c(\vect{\theta},\vect{\alpha},\mat{T})=\prod_{i=1}^\infty
\pi_i(\vect{\theta})^{L^i} \prod_{k=1}^p \alpha_k^{B_k^i}\prod_{\stackrel{k,\ell =1}{\ell\neq k}}^p
\left(\frac{t_{k\ell}}{s_i}\right)^{N_{k\ell}^i}
\exp \left( -\frac{t_{k\ell}}{s_i}Z_k^i \right) \prod_{k=1}^p
\left(\frac{t_{k}}{s_i}\right)^{N_{k}^i}
\exp \left( -\frac{t_{k}}{s_i}Z_k^i \right)  \]
with corresponding log--likelihood 
\begin{eqnarray}
  \ell_c (\vect{\theta},\vect{\alpha},\mat{T})&=&\sum_{i=1}^\infty L^i
\log \pi_i (\vect{\theta}) + \sum_{i=1}^\infty\sum_{k=1}^p B_k^i \log
\alpha_k + \sum_{i=1}^\infty \sum_{\stackrel{k,\ell=1}{\ell \neq k}}^p N_{k\ell}^i \log \left(
  \frac{t_{k\ell}}{s_i}\right)\nonumber  \\
&&-\sum_{i=1}^\infty \sum_{\stackrel{k,\ell=1}{\ell \neq k}}^p \frac{t_{k\ell}}{s_i} Z_k^i +
\sum_{i=1}^\infty \sum_{k=1}^p N_k^i \log \left( \frac{t_k}{s_i}\right)
-\sum_{i=1}^\infty \sum_{k=1}^p \frac{t_k}{s_i}Z_k^i . \label{eq:form-of-log-likelihood}
\end{eqnarray}
In order to calculate the maximum likelihood estimator
$(\hat{\vect{\theta}},\hat{\vect{\alpha}},\hat{\mat{T}})$, which is
the point that maximizes the likelihood (or log--likelihood)
function, we calculate first order partial derivatives of $t_{k\ell}$,
$t_k$, while we use the Lagrange multiplier methods
for $\alpha_i$ and $\pi_i(\vect{\theta})$ due to the
constraints on them summing up to one. 

Consider the parameter $t_{k\ell}$,  $k\neq \ell$. Then 
\[  \frac{\partial \ell_c}{\partial t_{k\ell}} =\sum_{i=1}^\infty
N_{k\ell}^i \frac{1}{t_{k\ell}}
-\sum_{i=1}^\infty  \frac{Z_k^i}{s_i}  =  0 \] 
implying
\[  \hat{t}_{k\ell} = \frac{\sum_{i=1}^\infty N_{k\ell}^i}{\sum_{i=1}^\infty Z_k^i/s_i}  .\]
Similarly, 
\[  \hat{t}_k = \frac{\sum_{i=1}^\infty N_k^i}{\sum_{i=1}^\infty Z_k^i/s_i}  . \]
The diagonal terms $\hat{t}_{kk}=-\sum_{\ell\neq k}\hat{t}_{k\ell} - \hat{t}_k$.
Regarding $\vect{\alpha}$, consider the Lagrange function 
\[ M(\vect{\alpha})=\sum_{i=1}^\infty B_k^i\log \alpha_k + \mu
(1-\sum_k \alpha_k) , \]
where $\mu$ is a Lagrange multiplier. Then 
\[ \frac{\partial M}{\partial \alpha_k} = \sum_{i=1}^\infty
\frac{B_k^i}{\alpha_k}-\mu =0 , \]
which result in
\[ \alpha_k \mu = \sum_{i=1}^\infty B_k^i .\]
Summing over $k$ yields
\[ \mu = \sum_{i=1}^\infty \sum_{k=1}^p B_k^i = M \]
so
\[ \hat{\alpha}_k = \frac{1}{M}\sum_{i=1}^\infty B_k^i  . \]
Concerning $\vect{\pi}(\vect{\theta})$, it depends on the particular
form of the discrete distribution whether it can be estimated
explicitly or numerically. We shall consider some particular examples.

We now consider the case of incomplete data. We shall employ the
EM--algorithm, which optimizes the incomplete data likelihood $L$,
\[  L(\vect{\theta},\vect{\alpha},\mat{T}; \vect{y})=\prod_{k=1}^M
f_Y(y_k;\vect{\theta},\vect{\alpha},\mat{T}) ,  \]
using the complete likelihood $L_c$ (or $\ell_c$) in the following
way. Here $f_Y(y;\vect{\theta},\vect{\alpha},\mat{T})$ denotes the density function of 
$Y\sim \mbox{NPH}_p(\pi(\vect{\theta}),\vect{\alpha},\mat{T})$.
  \begin{description}
  \item[0] Initialize with some ``arbitrary'' $(\vect{\theta}_0,\vect{\alpha}_0,\mat{T}_0)$ and let $n=0$.
  \item[1] (E--step) Calculate the function
    \[ h:(\vect{\theta},\vect{\alpha},\mat{T}) \rightarrow \Exp_{(\vect{\theta}_n,\vect{\alpha}_n,\mat{T}_n)}\left( \ell_c
      (\vect{\theta},\vect{\alpha},\mat{T}) | \vect{Y}=\vect{y} \right) .  \]
  \item[2] (M--step) Let $(\vect{\theta}_{n+1},\vect{\alpha}_{n+1},\mat{T}_{n+1}):=\mbox{argmax}_{(\vect{\theta},\vect{\alpha},\mat{T})}
    h(\vect{\theta},\vect{\alpha},\mat{T})$.
  \item[3] n=n+1; GOTO 1.
  \end{description}
In each iteration, the incomplete likelihood is increased (i.e. $L(\vect{\theta}_{n},\vect{\alpha}_{n},\mat{T}_{n})\leq L(\vect{\theta}_{n+1},\vect{\alpha}_{n+1},\mat{T}_{n+1})$)
and the procedure hence converges (possibly to a local maximum or saddlepoint though). 

We
notice that the actual calculations in the EM--algorithm only involve the
complete data likelihood. From the actual form of the log--likelihood
\eqref{eq:form-of-log-likelihood} we see that it is a linear function
of the sufficient statistics, and the conditional expected value of the
log--likelihood given the data will then be the log--likelihood function
with the sufficient statistics replaced by their conditional
expectations given the data. Hence the $M$--step is trivial since we only
have to plug in the conditional expectations given the data  instead of the sufficient
statistics which are not available

Concerning the E--step we proceed as follows.
First we consider one (generic) data point ($M=1$) and let $y=y_1$. We need to calculate the conditional expected
values of the sufficient statistics given $Y=y$. All distributions and expectations are under 
$\Prob=\Prob_{\vect{\Psi}}$, $\vect{\Psi}=(\vect{\theta},\vect{\alpha},\mat{T})$, but we will omit the index in order to ease the exposition. 

Concerning
$L^i$, it equals one if $i$ is the chosen level and zero otherwise. 
Let $I$ denote the
random variable indicating the chosen level. Then
\begin{eqnarray}
  \Exp (L^i | Y=y) &=& \Prob (I=i | Y=y) \nonumber \\
&=&\frac{\Prob (Y\in dy| I=i)\Prob (I=i)}{\Prob (Y\in dy)} \nonumber \\
&=& \frac{\pi_i (\vect{\theta}) \vect{\alpha}\exp (\mat{T}_i
  y)\vect{t}_i}{f_Y(y;\vect{\theta},\vect{\alpha},\mat{T})}  . \label{eq:orig-L}
\end{eqnarray}
Similarly, for $B_k^i$ we have
\begin{eqnarray*}
  \Exp (B_k^i | Y=y)&=&\Exp (1\{J_0=k, I=i\}| Y=y) \\
&=& \frac{\Prob (Y\in dy | J_0=k, I=i)\Prob (I=i)\Prob (J_0=k)}{\Prob
  (Y\in dy)} \\
&=& \frac{\pi_i (\vect{\theta}) \alpha_k \vect{e}_k^\prime \exp
  (\mat{T}_i y)\vect{t}_i}{f_Y(y;\vect{\theta},\vect{\alpha},\mat{T})} 
\end{eqnarray*}
so
\[  \Exp \left(\left.  \sum_{i=1}^\infty B_k^i \right| Y=y \right)=
\frac{\sum_{i=1}^\infty \pi_i (\vect{\theta}) \alpha_k \vect{e}_k^\prime \exp
  (\mat{T}_i y)\vect{t}_i}{f_Y(y;\vect{\theta},\vect{\alpha},\mat{T})} . \]
Regarding $Z_k^i$,
\begin{eqnarray*}
  \Exp \left( Z_k^i|Y=y \right)&=& \Exp \left( \left. \int_0^y 1\{
      J_u=k,I=i\} du\right| Y=y\right) \\
&=& \Exp \left( 1\{ I=i\} \left. \int_0^y 1\{
      J_u=k\} du\right| Y=y\right) \\
&=& \Exp \left( \left. \int_0^y 1\{
      J_u=k\} du\right| Y=y, I=i\right)\Prob (I=i|Y=y)\\
&=&\int_0^y \Prob (J_u=k|Y=y, I=i)du \frac{\Prob (Y\in dy|I=i)\Prob
  (I=i)}{\Prob (Y\in dy) } \\
&=&\int_0^y \frac{\Prob (Y\in dy, J_u=k|I=i)}{\Prob (Y\in dy|I=i)}du  \frac{\Prob (Y\in dy|I=i)\Prob
  (I=i)}{\Prob (Y\in dy) }\\
&=&\frac{\Prob (I=i)}{\Prob (Y\in  dy)} 
  \int_0^y {\Prob (Y\in dy, J_u=k|I=i)}du \\
&=&\frac{\Prob (I=i)}{\Prob (Y\in  dy)}
   \int_0^y {\Prob (Y\in dy| J_u=k, I=i) \Prob (J_u=k|I=i)}du \\
&=&\frac{\pi_i(\vect{\theta})}{f_Y(y;\vect{\theta},\vect{\alpha},\mat{T})}
  { \int_0^y\vect{e}_k^\prime e^{\mat{T}_i
    (y-u)}\vect{t}_i \vect{\alpha}
  e^{\mat{T}_i u}\vect{e}_k du } .
\end{eqnarray*}
The formula for $\Exp (N_{k\ell}^i|Y=y)$ is derived in a similar way,
resulting in
\[ \Exp (N_{k\ell}^i|Y=y) = \frac{
 t_{k\ell}}{s_i}\frac{\pi_i(\vect{\theta})
 }{f_Y(y;\vect{\theta},\vect{\alpha},\mat{T})}
\int_0^y \vect{e}_\ell^\prime e^{\mat{T}_i (y-u)}\vect{t}_i \vect{\alpha}
e^{\mat{T}_i u}\vect{e}_k du . \] 
Finally,
\[  \Exp (N_{k}^i|Y=y) = \frac{t_k}{s_i}
\frac{\pi_i(\vect{\theta})}{f_Y(y;\vect{\theta},\vect{\alpha},\mat{T})}
\vect{\alpha}e^{\mat{T}_iy}\vect{e}_k  . \]
The $\vect{\theta}$ needs to be treated separately on a case--by--case basis and it generally involves a numerical procedure for finding the
next iteration $\vect{\theta}^{(n+1)}$. We need to maximize
\[  \vect{\theta}\rightarrow \sum_{i=1}^\infty \Exp (L^i|Y=y) \log \pi_i
(\vect{\theta}) \]
subject to $\sum_{i=1}^\infty \pi_i(\vect{\theta})=1$, where (see \eqref{eq:orig-L})
\[ \Exp (L^i|Y=y)  = \frac{\pi_i (\vect{\theta})}{
{f_Y(y;\vect{\theta},\vect{\alpha},\mat{T})}}  
{ \vect{\alpha}
\exp \left(\mat{T}_i y\right)\vect{t}_i}  . \]
In general, for $M>1$ datapoints we simply sum the previous formulas with arguments $y_j$, $j=1,...,M$. 

\begin{Rem}\rm
We see that the formulas in the $E$--step involve both matrix--exponentials and integrals thereof, and by defining the generic integral
\[ \mat{J}(y):=\mat{J}(y;\vect{\alpha},\mat{T})=\int_0^y e^{\mat{T} (y-u)}\vect{t} \vect{\alpha}
e^{\mat{T} u} du ,  \]
we have
that (see \cite{VanLoan:1978tq}) 
\begin{equation}
 \exp \left( \matrix{\mat{T} & \vect{t}\vect{\alpha} \\ \mat{0} & \mat{T} } y \right) 
= \matrix{e^{\mat{T}y} & \mat{J}(y) \\ \mat{0} & e^{\mat{T}y} }  . \label{eq:non-unifomization} 
\end{equation}
Thus a simple (and numerically efficient) way of obtaining both $\exp (\mat{T}y)$ and $\mat{J}(y)$ is by calculating the matrix exponential on the left
hand side. 
\end{Rem}


The EM--algorithm can now be stated as follows.
\begin{Th}[EM--algorithm]\label{th:EM}
\mbox{ }
 \begin{description}
  \item[0] Initialize with some ``arbitrary'' $(\vect{\theta},\vect{\alpha},\mat{T})$.
  \item[1] (E--step) Compute 
\begin{eqnarray*}
 \Exp (L^i | \vect{Y}=\vect{y})
&=&\pi_i (\vect{\theta}) \sum_{j=1}^M \frac{\vect{\alpha}\exp (\mat{T}_i
  y_j)\vect{t}_i}{f_Y(y_j;\vect{\theta},\vect{\alpha},\mat{T})} \\
   \Exp (B_k^i | \vect{Y}=\vect{y})&=&\pi_i (\vect{\theta}) \sum_{j=1}^M
\frac{ \alpha_k \vect{e}_k^\prime \exp
  (\mat{T}_i y_j)\vect{t}_i}{f_Y(y_j;\vect{\theta},\vect{\alpha},\mat{T})}
\\
\Exp \left( Z_k^i|\vect{Y}=\vect{y} \right) &=&\pi_i(\vect{\theta})\sum_{j=1}^M
\frac{\mat{J}(y_j/s_i;\vect{\alpha}, \mat{T})_{kk}}{f_Y(y_j;\vect{\theta},\vect{\alpha},\mat{T})}
  \\
\Exp (N_{k\ell}^i|\vect{Y}=\vect{y}) &=& \frac{\pi_i(\vect{\theta})}{s_i}\sum_{j=1}^M
\frac{\mat{J}(y_j/s_i;\vect{\alpha}, \mat{T})_{\ell k}t_{k\ell} }{f_Y(y_j;\vect{\theta},\vect{\alpha},\mat{T})}
\\
\Exp (N_{k}^i|\vect{Y}=\vect{y}) &=&\frac{\pi_i(\vect{\theta})}{s_i} \sum_{j=1}^M 
\frac{\vect{\alpha}e^{\mat{T}_i y_j}\vect{e}_k t_k }{f_Y(y_j;\vect{\theta},\vect{\alpha},\mat{T})}
\end{eqnarray*}
  \item[2] (M--step) Maximize 
\[ \vect{\theta}\rightarrow \sum_{i=1}^\infty \Exp (L^i|\vect{Y}=\vect{y}) \log \pi_i
(\vect{\theta}) \]
subject to $\sum_i \pi_i(\vect{\theta})=1$ and let
$\hat{\vect{\theta}}$ be the argument at which the maximum is attained. 
Let 
\begin{eqnarray*}
  \hat{\alpha}_k&=&\frac{1}{M}\sum_{i=1}^\infty
  \Exp(B_k^i|\vect{Y}=\vect{y}) \\
  \hat{t}_{k\ell} &=& \frac{\sum_{i=1}^\infty \Exp
    (N_{k\ell}^i|\vect{Y}=\vect{y})}{\sum_{i=1}^\infty \frac{1}{s_i}\Exp (Z_k^i|\vect{Y}=\vect{y})}\\
 \hat{t}_k &=& \frac{\sum_{i=1}^\infty \Exp
   (N_k^i|\vect{Y}=\vect{y})}{\sum_{i=1}^\infty \frac{1}{s_i} \Exp (Z_k^i|\vect{Y}=\vect{y})} 
\end{eqnarray*}
and assign diagonal values $\hat{t}_{kk}=-\sum_{l\neq k}
\hat{t}_{k\ell}-\hat{t}_k$ so $\hat{\mat{T}}=\{ \hat{t}_{k\ell}\}_{k,\ell=1,...,p}$.
\item[3] Reassign values to initial parameters
  \begin{eqnarray*}
    \vect{\theta}&:=&\hat{\vect{\theta}} \\ 
    \vect{\alpha}&:=& \hat{\vect{\alpha}} \\
    \mat{T}&:=& \hat{\mat{T}} 
  \end{eqnarray*}
  \item[4] GOTO 1.
  \end{description}
\end{Th}

\begin{Rem}\label{rm:rem-repeated}\rm
Consider the $E$--step of Theorem \ref{th:EM}. Suppose that there are repeated values among the data points such that among the data $y_1,...,y_M$ there are $\tilde{y}_1,...,\tilde{y}_D$ different values and that $\tilde{y}_i$ appears $k_i\geq 1$ times in the original data. Then $k_1+\cdots k_D=M$ and the sums over $j=1,...,M$ in expected values can then be reduced to weighted sums of fewer terms instead. For example, 
\[ \Exp \left( Z_k^i|\vect{Y}=\vect{y} \right) =\pi_i(\vect{\theta})\sum_{j=1}^D k_j
\frac{\mat{J}(\tilde{y}_j/s_i;\vect{\alpha}, \mat{T})_{kk}}{f_Y(\tilde{y}_j;\vect{\theta},\vect{\alpha},\mat{T})} . \]
This shows that the EM--algorithm of Theorem \ref{th:EM} can also be used to estimate an NPH distribution when data are represented by a histogram rather than the raw data. This is particularly suitable when the amount of data is large since the histogram may here adequately represent the target distribution. 
\end{Rem}

\begin{Rem}\label{rm:rem-hist}\rm
In order to reduce the computational burden, we may replace the raw data by a histogram in certain regions of the support where data gather rather densely.  Heavy--tailed data with support on $[0,\infty )$ will typically concentrate in an interval $[0,T)$ (the "body" of the distribution) while it will be more scarce in $[T,\infty )$ (the "tail"). 

Now assume that $M$ is large and that the support for the data may be split into $[0,T)$ and $[T,\infty )$ for some $T>0$ such that the concentration of data points in $[0,T]$ is large. Then divide $[0,T]$ into $K$ disjoint sub--intervals 
$[T_i,T_{i+1})$, $i=0,1,\dotsc,K-1$  where $0=T_0<T_1<\cdots < T_K=T$ and $k_i$ denote the number of data points falling into $[T_i,T_{i+1})$. Then we may use the repeated data reduction of Remark \ref{rm:rem-repeated} for the interval $[0,T)$ by treating $\frac{T_i+T_{i+1}}{2}$, $i=0,...,K-1$ as data points with counts $k_i$. In fact we might choose any point in $[T_i,T_{i+1})$ as a representative, including any data point falling in this interval, however, if we choose the left points of the intervals, $T_i$, as representatives, we have to make 
special arrangements regarding the first interval in order not to provoke an atom at zero. 

Concerning the tail, $[T,\infty )$, we continue to use the raw data since the scattering in the tail is typically too diffuse in order to be represented by a histogram without too much loss of information.
\end{Rem}

\begin{Ex}\label{ex:explicit}\rm
Here we present an important special case for the choice of the scaling distribution $\{ \pi_i(\vect{\theta})\}$ of $N$ for the
Regularly Varying case.  It features an explicit solution to the M--step maximization problem of the function
\[  \theta\rightarrow \sum_{i=1}^\infty \log \pi_i(\theta)w_i  , \]
where $w_i=\Exp (L_i|\vect{Y}=\vect{y})$. Define
\begin{align*}
  \pi_i(\theta)=\Prob(N=s_i)=F(s_{i+1})-F(s_{i})&=\frac{1}{s_{i}^\theta}-\frac{1}{s_{i+1}^\theta},\qquad
  s_i=e^{(i-1)c}.\\
    &=(e^{-\theta c})^{(i-1)}(1-e^{-\theta c}),\qquad i=1,2,\dots
 \end{align*}
 This distribution corresponds to a the discretization of
 a Pareto distribution $F$ and the resulting discrete distribution is supported by $\{ e^{ic}: i=\in \mathds{N} \}$. For $c>0$, N with parameter 
 $e^{-\theta c}$.  The advantage of representing the scaling Pareto distribution on this form is that its argument which maximizes 
\[  \theta\rightarrow \sum_{i=1}^\infty \log \pi_i(\theta)w_i  , \]
where $w_i=\Exp (L_i|\vect{Y}=\vect{y})$, and it is not difficult to see it has an explicit solution given by
 \begin{equation*}
  \widehat{\theta}=-\dfrac 1c\,\log \left(1-\dfrac{\sum_{i=1}^\infty w_i}{\sum_{i=1}^\infty i\,w_i}\right) .
 \end{equation*}
The selected support of $N$ is also convenient for implementation purposes.
In practice, the infinite series defining the estimators can only be computed up to a finite 
number of terms.  The sequence $\pi(\theta)$ converges faster to $0$ so a significant reduction
{of terms to be computed is required in order to attain a certain specified precision}.

Since the distance between consecutive points in the support of $N$ is unbounded,
then the distribution of $N$ is not regularly varying and Breiman's lemma no longer applies.  
Nevertheless, the tail probability  of $Y$ oscillates between two regularly varying functions, so 
this model
can provide  an accurate approximation to any regularly varying distribution \cite{RojasXie}.
\end{Ex}

\section{Censoring}\label{sec:censoring}
In certain situations some data may be censored. 
Again we consider first the case of a single data point $y$ taken from a realization of $Y\sim \mbox{NPH}_p(\pi (\vect{\theta}),\vect{\alpha},\vect{T})$. The data point is right censored at $t$ if the only knowledge about $Y$ is that $Y>t$, left censored at $t$ if $Y\leq t$ and interval censored at $(s,t]$ if 
$Y\in (s,t]$. Left censoring is a special case of interval censoring with $s=0$ while right censoring can be obtained by fixing $s$ and letting $t\rightarrow \infty$. Formulas for right censoring will, however, appear as a part of the derivation of interval censoring.

The EM algorithm works entirely the same way as for uncensored data with the only difference that we are no longer observing a data point $Y=y$ but $Y\in (s,t]$. 
This will only change the $E$--step where we now have to calculate the following conditional expectations:
\begin{eqnarray*}
&& \Exp (L^i|Y\in (s,t]),\  \Exp (B_k^i|Y\in (s,t]),\ \Exp (Z_k^i|Y\in (s,t]), \ \Exp (N_{k\ell}^i | Y\in (s,t])\ \ \mbox{and} \ \ \Exp (N_k^i | Y\in (s,t]) .
\end{eqnarray*}
Concerning $L^i$, notice that
\begin{eqnarray*}
\Exp (L^i1\{ Y\leq t \})&=& \Prob (Y\leq t | I=i)\Prob (I=i) \\
&=&\pi_i(\vect{\theta})(1-\vect{\alpha}e^{\mat{T}_i t}\vect{e}) 
\end{eqnarray*}
and 
\[  \Exp (L^i1\{ Y>t \})= \pi_i(\vect{\theta})\vect{\alpha}e^{\mat{T}_it}\vect{e} .  \]
Thus 
\begin{eqnarray*}
\Exp (L^i|Y\in (s,t])&=& \frac{\Exp(L^i 1\{ Y \in (s,t]\})}{\Prob (Y\in (s,t])} \\
&=& \frac{\Exp (L^i) -\Exp(L^i 1\{ Y\leq s\}) - \Exp (L^i 1\{ Y>t \})}{\Prob (Y\in (s,t])} \\
&=& \frac{\pi_i(\vect{\theta})-\pi_i(\vect{\theta})(1-\vect{\alpha}e^{\mat{T}_is}\vect{e})-\pi_i(\vect{\theta})\vect{\alpha}e^{\mat{T}_it}\vect{e}}{\Prob (Y\in (s,t])} \\
&=& \frac{\pi_i(\vect{\theta})(\vect{\alpha}e^{\mat{T}_i s}\vect{e}-\vect{\alpha}e^{\mat{T}_it}\vect{e})}{\sum_{i=1}^\infty \pi_i(\vect{\alpha})(\vect{\alpha}e^{\mat{T}_is}\vect{e}-\vect{\alpha}e^{\mat{T}_i t}\vect{e})} .
\end{eqnarray*}
For right censored data we get 
\[ \Exp (L^i|Y>t) = \frac{\pi_i(\vect{\theta})\vect{\alpha}e^{\mat{T}_it}\vect{e}}{\Prob (Y>t)} = \frac{\pi_i(\vect{\theta})\vect{\alpha}e^{\mat{T}_it}\vect{e}}{\sum_{i=1}^\infty \pi_i(\vect{\theta})\vect{\alpha}e^{\mat{T}_it}\vect{e} } .  \]

The rest of the formulas are derived as for censored phase--type distributions (see \cite{Olsson:1996ua}) conditionally on the level $L^i$, with parameters $\vect{\alpha},\mat{T}_i$, which happens with probability $\pi_i(\vect{\theta})$. Thus we get that 
\begin{align*}
\Exp (B_k^i|Y\in (s,t])&= \frac{\pi_i(\vect{\theta})\left( \alpha_k \vect{e}_k^\prime e^{\mat{T}_is}\vect{e}-\alpha_k \vect{e}_k^\prime e^{\mat{T}_ie}\vect{e}\right) }{\Prob (Y\in (s,t])} \\
\Exp (N_{k\ell}^i|Y\in (s,t])&= \dfrac{\pi_i(\vect{\theta})t_{k\ell}/s_i \left(\int_s^t\vect{\alpha}e^{\mat{T}_iu}\vect{e}_k du -\int_s^t \vect{e}_{\ell}^\prime e^{\mat{T}_i(t-u)}\vect{e}\vect{\alpha}e^{\mat{T}_iu}\vect{e}_kdu \right)}{\Prob (Y\in (s,t])} \\
\Exp (Z_{k}^i|Y\in (s,t])&= \dfrac{\pi_i(\vect{\theta}) \left(\int_s^t\vect{\alpha}e^{\mat{T}_iu}\vect{e}_k du -\int_s^t \vect{e}_{k}^\prime e^{\mat{T}_i(t-u)}\vect{e}\vect{\alpha}e^{\mat{T}_iu}\vect{e}_kdu \right)}{\Prob (Y\in (s,t])} \\
\Exp (N_{k}^i|Y\in (s,t])&= \dfrac{\pi_i (\vect{\theta}) t_k/s_i \int_s^t \vect{\alpha}e^{\mat{T}_iu}\vect{e}_kdu}{\Prob (Y\in (s,t])}
\end{align*} 
with similar and obvious formulas for the right censored case. The integrals are calculated similarly as in the uncensored case, namely
\begin{eqnarray*}
\int_s^t \vect{\alpha}e^{\mat{T}_iu}\vect{e_k}du
 &=& 
\vect{\alpha}\mat{T}_i^{-1}(e^{\mat{T}_is}-e^{\mat{T}_it})\vect{e}_k \\
&=&s_i \vect{\alpha}\mat{T}^{-1}(e^{\mat{T}_is}-e^{\mat{T}_it})\vect{e}_k \\
\int_s^t\vect{e}_{\ell}^\prime e^{\mat{T}_i(t-u)}\vect{e}\vect{\alpha}e^{\mat{T}_iu}\vect{e}_kdu&=&\mat{J}^i(t)_{\ell k}-\mat{J}^i(s)_{\ell k},
\end{eqnarray*}
where
\[  \mat{J}^i(y)=\exp \left[  
\begin{pmatrix}
\mat{T}_i & \vect{e}\vect{\alpha} \\
\mat{0} & \vect{T}_i 
\end{pmatrix}
y\right] .\]

For more than one data point, the data are split into a group of uncensored data and into other groups of different types of censored data. The  conditional expectations are then calculated for all data points subject to their group classification,
 and  all conditional expectations of the same kind (jumps, occupation times etc.) are then summed over all data. This amounts to the E--step in an EM algorithm, the rest of which is identical to Theorem \ref{th:EM}. 

\begin{Rem}\rm
In Remark \ref{rm:rem-hist} we suggested the use of a histogram for representing the ``body'' of the distribution if the amount of data is large. Interval censoring provides a feasible alternative in this same direction. 
\end{Rem}

\section{Fitting to a known distribution}\label{sec:KL}
It may be of interest to approximate a given (heavy tailed) distribution $H$ by a distribution 
$G\in \mbox{NPH}$ if for example methods for calculation the ruin probability based on the claim size distribution $H$ is not known. In \cite{ANO:1996} it was shown how the EM algorithm can be modified in order to approximate phase--type distributions to a given distribution $F$ with non--negative support. The EM algorithm then converges to a limit which minimizes the Kullback--Leibler distance between phase--type distributions of a specified order and the distribution $F$.

The idea is to let the number of data points $M\rightarrow \infty$ such that the distribution $H$ is seen as an empirical distribution function for a data set with $M=+\infty$. If we let $M\rightarrow \infty$ then (see Theorem \ref{th:EM})
\[  \hat{\alpha}_k = \sum_{i=1}^\infty \pi_i (\vect{\theta}) \frac{1}{M}\sum_{j=1}^M
\frac{ \alpha_k \vect{e}_k^\prime \exp
  (\mat{T}_i y_j)\vect{t}_i}{f_Y(y_j;\vect{\theta},\vect{\alpha},\mat{T})} \rightarrow \sum_{i=1}^\infty \pi_i (\vect{\theta}) \int_0^\infty \frac{ \alpha_k \vect{e}_k^\prime \exp
  (\mat{T}_i y)\vect{t}_i}{f_Y(y;\vect{\theta},\vect{\alpha},\mat{T})}dH(y) . \]
Similarly it is not difficult to see that
\begin{align*}
\hat{t}_{k\ell}&= \frac{\sum_{i=1}^\infty \Exp
    (N_{k\ell}^i|\vect{Y}=\vect{y})}{\sum_{i=1}^\infty \frac{1}{s_i}\Exp (Z_k^i|\vect{Y}=\vect{y})} 
   = \frac{\sum_{i=1}^\infty \frac{1}{M}\Exp
    (N_{k\ell}^i|\vect{Y}=\vect{y})}{\sum_{i=1}^\infty \frac{1}{s_i}\frac{1}{M}\Exp (Z_k^i|\vect{Y}=\vect{y})} \\
    &\rightarrow \dfrac{\sum_{i=1}^\infty \frac{\pi_i(\vect{\theta})}{s_i}\int_0^\infty
\frac{\mat{J}(y/s_i;\vect{\alpha}, \mat{T})_{\ell k}t_{k\ell} }{f_Y(y;\vect{\theta},\vect{\alpha},\mat{T})}dH(y)}{\sum_{i=1}^\infty \frac{\pi_i(\vect{\theta})}{s_i}\int_0^\infty
\frac{\mat{J}(y/s_i;\vect{\alpha}, \mat{T})_{kk} }{f_Y(y;\vect{\theta},\vect{\alpha},\mat{T})}dH(y)}
\end{align*}
and
\begin{align*}
\hat{t}_k&\rightarrow\dfrac{\sum_{i=1}^\infty \frac{\pi_i(\vect{\theta})}{s_i}\int_0^\infty
\frac{\vect{\alpha}e^{\mat{T}_i y}\vect{e}_k t_k}{f_Y(y;\vect{\theta},\vect{\alpha},\mat{T})}dH(y)
}{\sum_{i=1}^\infty \frac{\pi_i(\vect{\theta})}{s_i}\int_0^\infty
\frac{\mat{J}(y/s_i;\vect{\alpha}, \mat{T})_{kk} }{f_Y(y;\vect{\theta},\vect{\alpha},\mat{T})}dH(y)} .
\end{align*}
Concerning $\vect{\theta}$, maximizing 
\[ \vect{\theta}\rightarrow \sum_{i=1}^\infty \Exp (L^i|\vect{Y}=\vect{y}) \log \pi_i
(\vect{\theta})  \]
is equivalent to maximizing 
\[ \vect{\theta}\rightarrow \sum_{i=1}^\infty \frac{1}{M}\Exp (L^i|\vect{Y}=\vect{y}) \log \pi_i
(\vect{\theta}) .  \]
As $M\rightarrow \infty$ the latter converges to 
\[  \vect{\theta}\rightarrow \sum_{i=1}^\infty \left( \int_0^\infty \frac{\pi_i (\vect{\theta}) \vect{\alpha}\exp (\mat{T}_i
  y)\vect{t}_i}{f_Y(y;\vect{\theta},\vect{\alpha},\mat{T})}dH(y)\right) \log \pi_i
(\vect{\theta}) , \]
which is then the function to be maximized.
Let $\hat{\vect{\theta}}$ denote the argument which maximizes this function.

In general none of the integrals will have explicit solutions,  and approximations (e.g.\ numerical integration, Quasi Monte Carlo methods
 or more sophisticated variants) will have to be employed.
The EM algorithm now works as follows: 
\begin{description}
\item[0] initiate with some parameters $(\vect{\theta},\vect{\alpha},\mat{T})$.
\item[1] calculate $\hat{\vect{\alpha}}$, $\hat{\mat{T}}$ and $\hat{\vect{\theta}}$.
\item[2] assign $(\vect{\theta},\vect{\alpha},\mat{T})=(\hat{\vect{\theta}},\hat{\vect{\alpha}},\hat{\mat{T}})$. 
\item[3] GOTO 1.
\end{description}

\section{Examples}\label{sec:numerical}
{In this section we provide five examples. In the first example we consider the simplest case of simulated data from a scale mixture of exponential distributions, where the scaling distribution has a Pareto type of tail. In the second example we fit NPH distributions to a real data set (Danish reinsurance data of fire claims), while in the last three examples we consider the fitting of NPH distributions to the theoretical distributions log--gamma, Weibull and log--normal respectively.}


\begin{Ex}[Erlang distributions]\rm
A $q$ dimensional Erlang distribution, $\mbox{ER}_q(\lambda)$, is a phase--type distribution with canonical representation
\[ \vect{\alpha}=(1,0,0,...,0), \ \ \mat{T}=\left(
 \begin{array}{ccccc}
  -\lambda & \lambda & 0 & \cdots & 0 \\
    0 & -\lambda & \lambda & \cdots & 0 \\
    0 & 0 & -\lambda & \cdots & 0 \\
    \vdots & \vdots & \vdots & \vdots\vdots\vdots & \vdots \\
    0 & 0 & 0 & \cdots & -\lambda
 \end{array}
\right) \ \ \mbox{and}\ \ \vect{t} =\left(
 \begin{array}{c}
   0 \\
  0 \\
   0 \\
   \vdots \\
    \lambda
 \end{array}
\right)  .  \]
Hence, the correspoding NPH distribution has density
\begin{equation}\label{mixture.erlangs}
 f_Y(y;\vect{\theta},\lambda)=
\sum_{i=1}^\infty\pi_i(\vect{\theta})h(y/s_i;q,\lambda)
=\sum_{i=1}^\infty
\pi_i(\vect{\theta}) \left(\frac{\lambda}{s_i}\right)^q
\frac{y^{q-1}}{(q-1)!}e^{-y\lambda /s_i},
\end{equation}
and the maximum likelihood estimator for $\lambda$ is easily calculated to be
\[  \widehat{\lambda}^{-1}= \frac{1}{qM}\sum_{i=1}^\infty \frac{Z^i}{s_i},\qquad Z^i=\sum_{j=1}^N y_j\cdot 1(I_j=i),  \]
while the maximum likelihood estimate for $\vect{\theta}$ has the general form
\begin{equation*}
 \widehat{\vect{\theta}}=\argmax_{\theta}\sum_{i=1}^\infty L^i\log\pi_i({\theta}),\qquad  L^i=\sum_{j=1}^N 1(I_j=i).
\end{equation*}
The $E$--step is similar as for the general model.
In the case of the sufficient statistic $Z^i$ we have
\begin{align*}
 \Exp (Z^i|\vect{Y}=\vect{y})= \sum_{j=1}^M y_j\Prob (I_j=i|Y=y_j)
  &=\sum_{j=1}^M y_j \frac{\Prob (I_j=i)\Prob (Y\in dy_j|I=i)}
  {\Prob (Y\in dy_j)}\\
  &=\sum_{j=1}^M y_j \frac{\pi_i(\vect{\theta})
  h(y_j;q,\lambda /s_i)}
{f(y_j;\vect{\theta} ,\lambda)} .
\end{align*}

As for the sufficient statistic $L^i$, its expected value is analogue for the general model in the previous section and found to be equal to
  \begin{equation*}
    \Exp \left( L^i |\vect{Y}=\vect{y}\right)
    = \sum_{j=1}^M \frac{ \pi_i(\vect{\theta})
      h(y_j/s_i;q,\lambda)}
      {f(y_j;\vect{\theta},\lambda)}  .
  \end{equation*}

We present a small simulation study for the case of $q=1$ (corresponding to infinite dimensional hyperexponential distribution), $s_i=i$ and scaling distribution 
\[ \pi_i (\theta) = \frac{i^{-\theta}}{\zeta (\theta)}, i=1,2,... ,\]
where
  \[ \zeta (\theta ) = \sum_{i=1}^\infty i^{-\theta} \] 
  is the
  Riemann Zeta function with parameter $\theta$. This distribution is also known as the Riemann Zeta  or discrete Pareto distribution since its tail resemble that of a Pareto distribution.

  In order to find the EM--estimate of $\theta^{(n+1)}$, we have to maximize the function
  \[ h(\theta) = -\sum_{i=1}^\infty \Exp (L^i|\vect{Y}=\vect{y})
  \theta \log (i) - \sum_{i=1}^\infty \Exp (L^i|\vect{Y}=\vect{y})
  \log (\zeta (\theta )) .\] 
  Differentiating with respect to $\theta$
  and rearranging implies that we have to solve the following equation
  w.r.t. $\theta$,
  \begin{equation}
    \frac{\zeta^\prime (\theta )}{\zeta (\theta )} =
    -\frac{1}{M}\sum_{i=1}^\infty \Exp \left( L^i
      |\vect{Y}=\vect{y}\right) \log (i)  ,\label{eq:MLi}
  \end{equation}
  which is done numerically using a simple Newton--Raphson procedure.

  Results of the simulation study, which is shown in
  Table \ref{table:exponential}, reveals that the EM algorithm
  is able to recover the underlying structure of the data, and that
  the estimation, as expected, improves with larger sample sizes.
  \begin{table}[ht]
    \begin{center}
      \begin{tabular}{|cc|cc|cc|cc|cc|cc|}
        \multicolumn{2}{|c|}{parameters} & \multicolumn{10}{|c|}{sample size}\\ \hline
        $\lambda$ & $\theta$ &  \multicolumn{2}{|c|}{$100$} &
        \multicolumn{2}{|c|}{$500$} &  \multicolumn{2}{|c|}{$1000$} &
        \multicolumn{2}{|c|}{$5000$} &
        \multicolumn{2}{|c|}{$10000$}  \\ \hline
        1.0 & 2.0 & 1.04 & 1.88 & 0.95 & 2.09  & 1.14 & 1.93 & 0.99 & 2.00 & 1.01 & 2.01 \\
        1.0 & 2.5 & 0.85 & 2.67 & 0.90 & 2.52 & 1.01   & 2.63  & 0.94 & 2.64
        & 1.00 & 2.52  \\
        1.0 & 5.0 & 0.92 & 7.1 & 1.09 & 7.26  & 1.05 & 4.65 & 1.00 & 4.85 &
        1.01 & 5.03 \\
        \hline
      \end{tabular}
    \end{center}
    \caption{\label{table:exponential} EM--estimates of
      $(\widehat{\lambda},\widehat{\theta})$
      infinite
      dimensional hyper--exponential distributions with varying parameters
      and sample size}
  \end{table}
\end{Ex}

\begin{Ex}\rm
We consider 2167 reinsurance data for Danish fire insurance claims above 1 million DKR for the period 1980--1993. 
These data have been widely studied in Extreme Value Theory  \cite{McNeil:1997,Resnick:1997}.
 The data corresponds to claims in millions of Danish Kroner for the period 1980--1993, 
 the amounts being adjusted for inflation to prices of 1985. We subtracted 1 (million) from all data in order to shift
  the support to $[0,\infty)$ which is the natural support for a phase--type distributions. Of the 2167 data, 519 are repeated values so only 1648 have different values.

We propose an NPH model for the data 
by mixing a five dimensional phase--type distribution with
  an $N$ that is assumed to follow the discretized Pareto distribution of Example \ref{ex:explicit}.

Just above 90$\%$ of the data is below 5, while less than 10$\%$ falls between 5 and the maximum observation of 262.2504. Thus we select 
the ``body part'' of the distribution as the densely populated region  $[0,5]$
while the much more diffuse region $[5,\infty)$ is then considered the ``tail''. 

First we test the difference in performance between the EM--algorithm using raw data and one using the hybrid method proposed in Remark \ref{rm:rem-hist}. 

We divide $[0,5]$ into 100 subintervals of the same size, thus treating losses within 50,000 DKR as all being the same. We choose the centre  of intervals as data points. In the tail we have 186 data points of which 8 are repeated values so a total of 178 different values. This amounts to a total of 278 distinct data points. The result of the experiment is following. While the actual discretization of $[0,5]$ does not appear to have any effect at all on the estimation, as compared with the EM--algorithm of the raw data, the speed was increased by a factor 5.75 which is more or less consistent with the execution times depending linearly on the number of data points, in which case we should have expected an increase in speed by a factor $1648/278=5.92$. 

In all the EM algorithms which follows we have used the hybrid method of Remark \ref{rm:rem-hist} with $[0,5]$ 
divided into 100 equally sized subintervals.

We now investigate the effect the choice of the parameter $c$ (see Example \ref{ex:explicit}) will have on the estimation. Intuitively, the smaller we choose the discretization steps of the geometric progression, the better the estimation. 
 Fixing the dimension of the phase--type 
 distribution to five, we select three different values for $c$, $1$, $1/2$ and $1/4$.  
 To obtain the EM estimates we randomly selected various sets of initial and iterated the EM until the relative change in 
 the loglikelihood was smaller than $10^{-8}$. We kept the model with the largest likelihood. 
 The estimates are as follows 
\begin{description}
\item[$c=1$] 
\begin{eqnarray*}
\hat{\theta}&=&1.2743\\
\hat{\vect{\alpha}}&=&(0.6415,   0.0099,  0.0055,  0.2115 ,  0.1316)\\
\hat{\mat{T}}&=&
\begin{pmatrix}
   -2.7430 &   1.3565 &        0  &       0  &       0\\
    0.0003 &  -3.0398 &        0  &       0  &       0\\
         0 &        0 &  -2.6313  &  1.1226  &  0.2167\\
         0 &        0 &   0.7223  & -1.3953  &  0.4089\\
         0 &        0 &   1.5129  &  0.8388  & -2.4779
\end{pmatrix}
\end{eqnarray*}
\item[$c=1/2$] 
\begin{eqnarray*}
\hat{\theta}&=&1.3136\\
\hat{\vect{\alpha}}&=&( 0.7692,    0.0149,    0.1154,    0.0148,    0.0857)\\
\hat{\mat{T}}&=&
\begin{pmatrix}
   -3.0713 &   0.1234 &   0.2319 &   0.0655 &   1.6457\\
    0.2239 &  -1.9576 &   0.3787 &   0.9651 &   0.3899\\
    0.1824 &   0.7410 &  -3.0705 &   0.7978 &   0.4776\\
    0.3099 &   0.6150 &   0.4751 &  -1.5588 &   0.1588\\
    0.0034 &   0.0446 &   0.2263 &   0.0446 &  -3.4154
\end{pmatrix}
\end{eqnarray*}
\item[$c=1/4$] 
\begin{eqnarray*}
\hat{\theta}&=&1.3230\\
\hat{\vect{\alpha}}&=&(0.0267,    0.4563,    0.0010,    0.2603,    0.2556)\\
\hat{\mat{T}}&=&
\begin{pmatrix}
   -0.7896 &   0.2531 &   0.0385 &   0.2652 &   0.0265\\
    0.0810 &  -3.7022 &   1.7356 &   1.1626 &   0.3934\\
    0.5838 &   0.0013 &  -3.7890 &   0.0576 &   0.0302\\
    0.3386 &   0.0606 &   1.0320 &  -3.7735 &   0.1957\\
    0.6694 &   0.0292 &   0.2376 &   0.3920 &  -3.4129
\end{pmatrix}
\end{eqnarray*}
\end{description}
That the estimates for $(\vect{\alpha},\mat{T})$ look distinct is because phase--type representations are by no means unique, however, the estimated densities which are plotted in Figure \ref{Danish.Dens} against the histogram of the data 
are almost identical (actually indistinguishable in the plotted range of $[0,6]$). It is clear that the different phases do not have a physical interpretation but are merely dummy (or black box) states used for the only purpose of obtaining a proper adjustment of the distribution.

In order to inspect the tail behaviour we plotted the log--survival function in Figure \ref{Danish.Survival}. We have provided two plots of the survival function in different ranges in order to be able to inspect the tail behaviour well beyond the largest data value (which is 262.2504). 
The different values for $\theta$ (obtained as a consequence of the different choices of $c$) does only result in a slightly different tail behaviour for very large claim sizes. We conclude that the estimation does not depend significantly on the choice of $c$ within the possible choices of $c=1,1/2$ or $1/4$, and can probably be extrapolated to concluding robustness against choices in $c$ within any ``reasonable'' range for $c$.

\begin{figure}[ht]
  \includegraphics[trim={5cm 13cm 5cm 12.5cm},clip,width=0.80\textwidth]{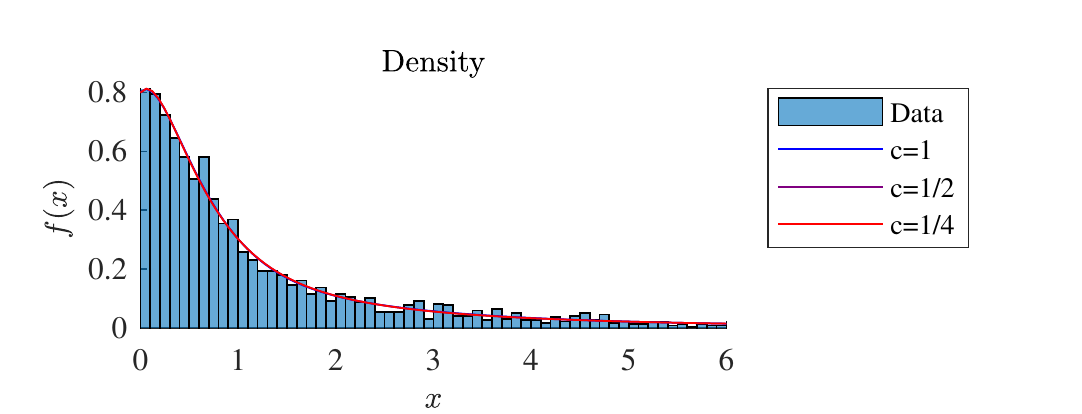}
  \caption{Estimated densities against histogram of the NPH models for Danish reinsurance data.}
  \label{Danish.Dens}
 \end{figure}

  \begin{figure}[ht]
  \includegraphics[trim={5cm 13cm 5cm 12.5cm},clip,width=0.80\textwidth]{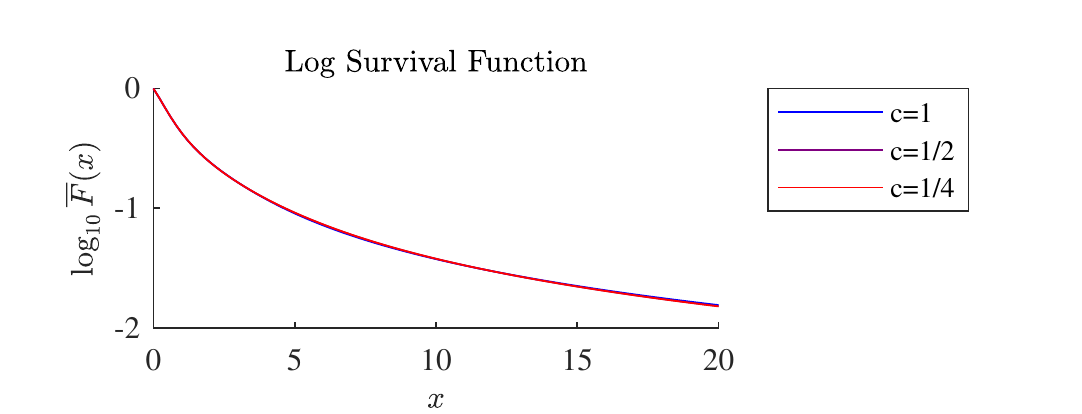}
  \includegraphics[trim={5cm 13cm 5cm 12.5cm},clip,width=0.80\textwidth]{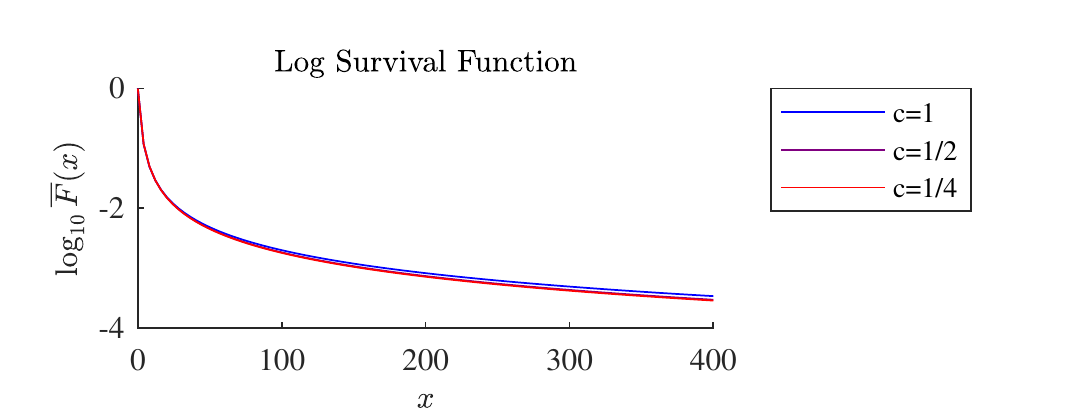}
  \caption{Estimated survival functions of the NPH models for Danish reinsurance data corresponding to $\theta$ equal to $1.2748, 1.3115, 1.3387$.}  \label{Danish.Survival}
 \end{figure} 
  Finally we compare the fit for $c=1/4$ to one obtained by \cite{McNeil:1997} and \cite{Resnick:1997}. 
 The first reference fits a Generalized Pareto Distribution via maximum likelihood estimation while the second employs the Hill estimator.
 The implementation of the methods described above involve the selection of an appropriate threshold
 which should be chosen by the modeler.  According to \cite{Resnick:1997} the values 
 in the interval [1.40, 1.46] produce excellent fits, the value recommend by the same author being $1.45$.

 We finally compare our estimated model with $c=1/4$ to a NPH model where we fix $\theta=1.45$ as recommended by 
  \cite{Resnick:1997}. This can be done by running the EM algorithm in the usual way but avoiding any adjustments in
  $\theta$. The results of the estimation are given next
\begin{description}
\item[$\theta=1.45$] Fixed tail. 
\begin{eqnarray*}
\hat{\vect{\alpha}}&=&( 0.1114,    0.0080,    0.3366,    0.3600,    0.1840)\\
\hat{\mat{T}}&=&
\begin{pmatrix}
   -3.0292 &   0.8810 &   0.0655 &   0.1600 &   0.1538\\
    0.1204 &  -0.6765 &   0.3295 &   0.0829 &   0.0827\\
    1.4308 &   0.2854 &  -4.2611 &   1.1133 &   1.4316\\
    0.4463 &   0.2717 &   0.3082 &  -3.6472 &   1.3508\\
    0.0958 &   0.0739 &   0.0230 &   0.0265 &  -3.2740
\end{pmatrix}
\end{eqnarray*}
\end{description}
Again this representation does not resemble any of the previous estimations (due to non--uniqueness of representations) but from figures \ref{Danish.Dens1} and \ref{Danish.Survival1} it is clear that there is no distinguishable difference between densities for the fixed and EM adjusted tails in the range $[0,6]$ where the main body of the distribution is situated. The tail behaviour also looks quite similar though the EM fitted 
tail seems to be slight heavier than the fixed tail.

 
 \begin{figure}[ht]
  \includegraphics[trim={5cm 13cm 5cm 12.5cm},clip,width=0.80\textwidth]{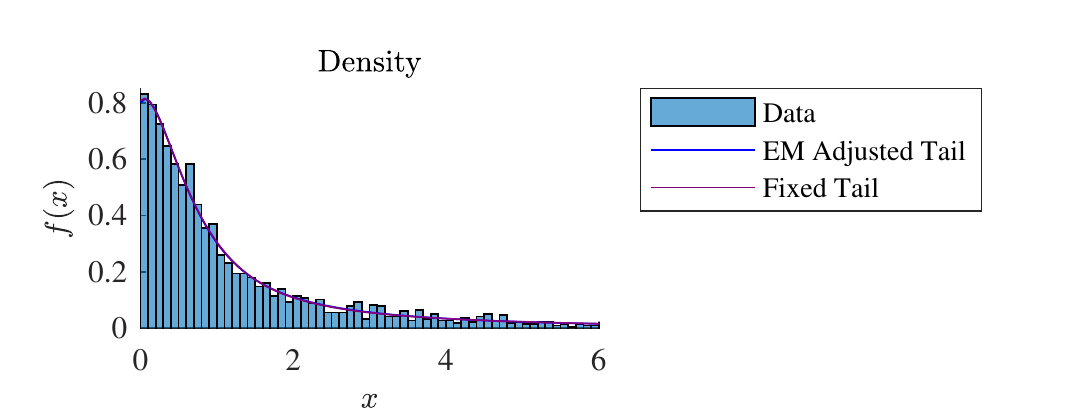}
  \caption{Estimated densities of the NPH models with EM adjusted ($\hat{\theta}=1.3387$) and fixed ($\hat{\theta}=1.45$) tails against 
  histogram of the Danish reinsurance data.  
  }
  \label{Danish.Dens1}
 \end{figure}

  \begin{figure}[ht]
  \includegraphics[trim={5cm 13cm 5cm 12.5cm},clip,width=0.70\textwidth]{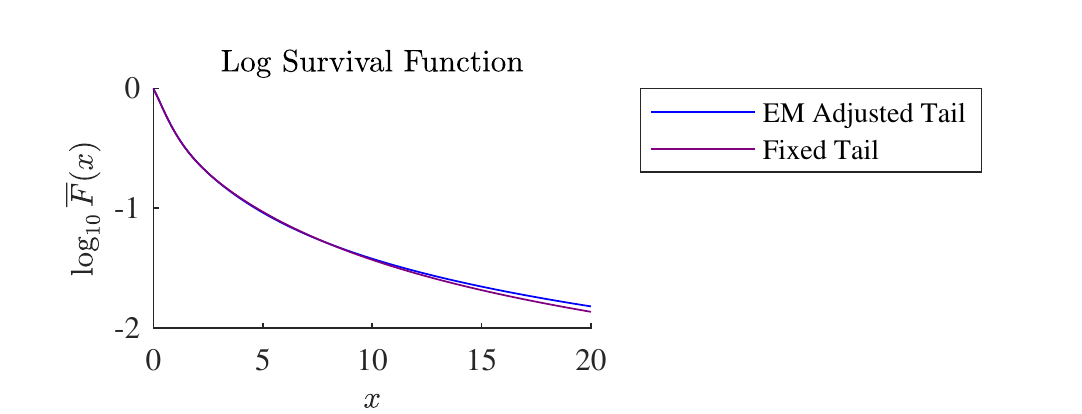}
  \includegraphics[trim={5cm 13cm 5cm 12.5cm},clip,width=0.70\textwidth]{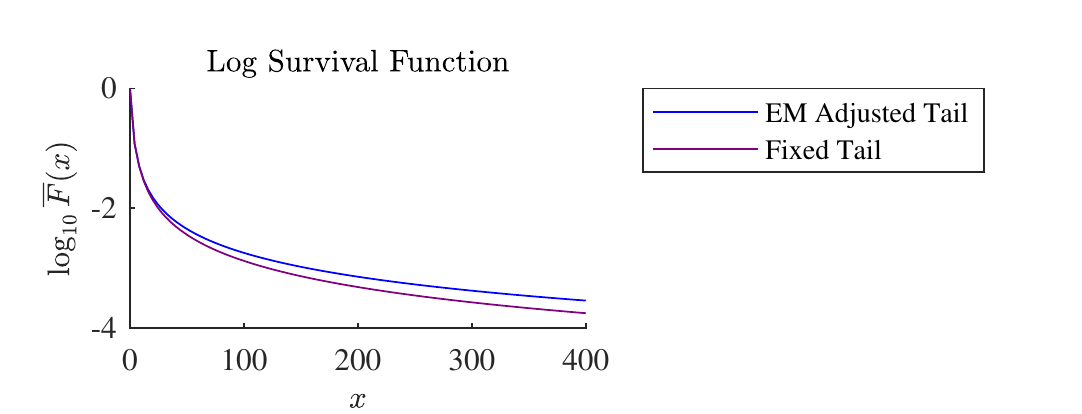}
  \caption{Estimated survival functions of the NPH models with EM adjusted ($\hat{\theta}=1.3387$) and fixed ($\hat{\theta}=1.45$) tails.}  \label{Danish.Survival1}
 \end{figure}   
\end{Ex}

\begin{Ex}[Fitting to a theoretical distribution: Loggamma]\rm
Next we consider the problem of approximating a theoretical distribution via an NPH model.
In our first example we take as target a log--gamma distribution
with shape parameter $\alpha$, scale parameter $\beta$  and shifted
one unit to the left so its support is $[0,\infty)$, which is the natural support for the class NPH. Its density function is then given by
\begin{equation*}
 f(x)=\dfrac{\beta^\alpha\log ^{\alpha-1}(x+1)}{\Gamma(\alpha)(x+1)^{\beta+1}},\qquad x\ge0.
\end{equation*}   
This distribution is regularly varying with parameter $\beta$. For this example we choose the parameters $\alpha=\beta=2$, 
for the target distribution with the purposes of analysing 
a distribution with a mode away from $0$ and a moderately heavy--tail. 
We consider an NPH model where the underlying  phase--type distribution has five phases
and the scaling distribution is that of Example \ref{ex:explicit} with $c=1$. With this 
example, we want to test if general Regularly Varying distributions can be correctly fitted
with the general model suggested in Example \ref{ex:explicit}.

We employed Quasi Monte Carlo ideas to approximate the integrals in the EM steps and iterated
the algorithm until the relative error was smaller than $10^{-9}$.
The results are given below
\begin{eqnarray*}
\hat{\theta}&=&1.6031\\
\hat{\vect{\alpha}}&=&(  0.5717,    0.0330,    0.0000,    0.3954,    0.0000)\\
\hat{\mat{T}}&=&
\begin{pmatrix}
   -1.9634 &   0.0609 &   0.5025 &   0.1249 &   1.2751\\
    0.0616 &  -0.3372 &   0.0775 &   0.0382 &   0.1428\\
    0.7529 &   0.1178 &  -2.2723 &   0.4797 &   0.0068\\
    0.7278 &   0.3060 &   1.1458 &  -4.8966 &   2.7170\\
    0.8923 &   0.0317 &   0.0482 &   0.2021 &  -3.4321
\end{pmatrix}
\end{eqnarray*}  

The densities of the log-gamma distribution and its NPH approximation are plotted in Figure \ref{Loggamma.Dens}.
The densities of the log-gamma and its NPH approximation are almost indistinguishable from each other in the body region.
The tail behavior is correctly captured by the NPH model as seen in Figure \ref{Loggamma.Survival1}.
The shape of the tail of the NPH is very close to that of the Loggamma distribution although the NPH estimate has
a heavier tail (the estimated regularly varying parameter was $\hat\theta=1.6031$).

 \begin{figure}[ht]
  \includegraphics[trim={5cm 13cm 5cm 12.5cm},clip,width=0.80\textwidth]{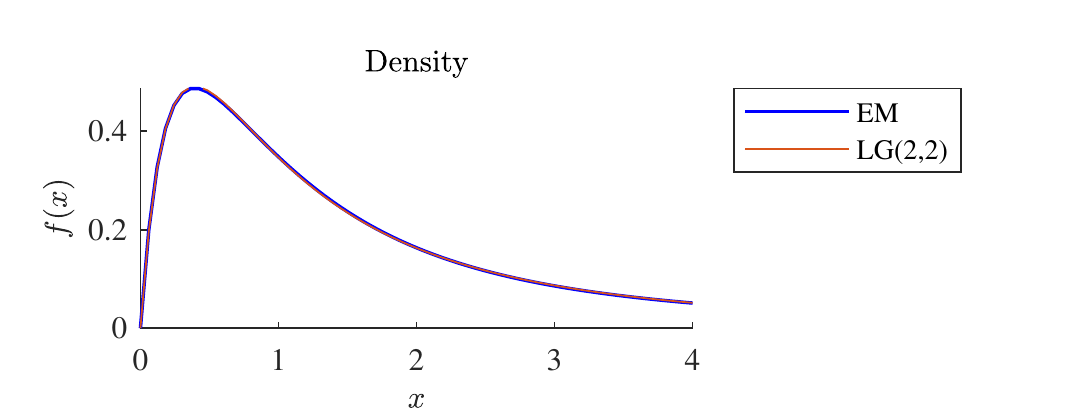}
  \caption{Density of the Loggamma distribution with parameters $(2,2)$ and an approximation via 
   an NPH model with EM adjustment.  
  }
  \label{Loggamma.Dens}
 \end{figure} 

  \begin{figure}[ht]
  \includegraphics[trim={5cm 13cm 5cm 12.5cm},clip,width=0.70\textwidth]{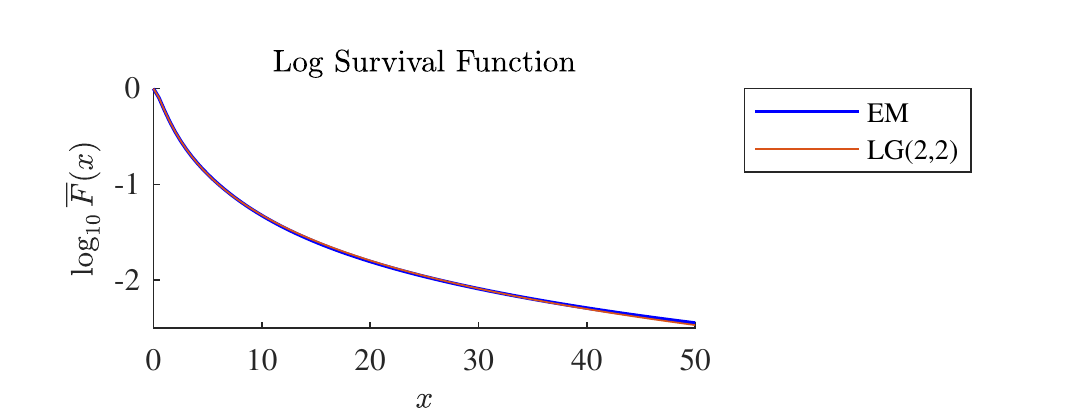}
  \includegraphics[trim={5cm 13cm 5cm 12.5cm},clip,width=0.70\textwidth]{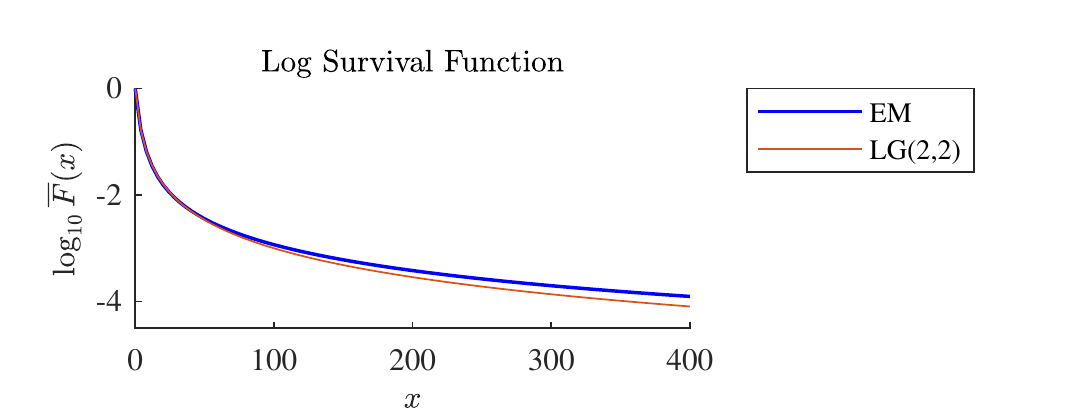}
  \caption{Survival functions of the Loggamma distribution 
   $\mbox{LG}(2,2)$ and an approximation via an NPH model with EM adjustment.}  \label{Loggamma.Survival1}
 \end{figure}

\end{Ex}

\begin{Ex}[Fitting to a theoretical distribution: Weibull]\rm
Next we move away from the Regularly Varying case and consider instead the Weibullian case.
As a target distribution we consider a classical two--parameter
Weibull with $\lambda=1$ and $p=1/2$ so its density is
$e^{-\sqrt{x}}/2\sqrt{x}$, $x\ge0$.
For adjusting this model we consider an NPH family of distributions where the phase--type part 
has  five phases and the scaling distribution is supported over a geometric progression 
$e^{c},e^{2c},e^{3c},\dots$ with $c=1$ and taken as 
a discretization of a classical two--parameter Weibull distribution with $\delta=p-1$
More precisely, its density is given by
\begin{align*}
 f(x)=p\lambda^p x^{p-1} e^{-(\lambda x)^p},\qquad x>0.
\end{align*}
Notice, that the target distribution is in the same two--parameter family of Weibull
distributions. The results are given below.
\begin{eqnarray*}
\hat{\lambda}&=&0.6181,\ \ \hat{p}=1.1673\\
\hat{\vect{\alpha}}&=&(  0.2370,    0.3349,    0.0901,    0.3380,         0)\\
\hat{\mat{T}}&=&
\begin{pmatrix}
   -3.0858 &   0.6964 &   0.2188 &   0.3355 &   0.5572\\
   51.0253 &-207.2799 &  18.3961 &  18.5851 &  46.0724\\
    0.4369 &   0.2511 &  -0.9839 &   0.0122 &   0.0487\\
    1.1297 &   0.2448 &   0.9388 & -11.0548 &   0.7637\\
    0.7583 &   1.0994 &   0.9844 &   0.4882 &  -3.3303
\end{pmatrix}
\end{eqnarray*}  
The agreement between the Weibull distribution and its NPH approximation is very good in the body of the distribution.
The approximation of the tail is also particularly good for values going up to 100 which
correspond to probabilities of order $10^{-4}$. The NPH adjusted model is Weibullian with 
parameter $\widehat p(1+\widehat p)^{-1}\approx 1.1673/2.1673\approx0.5386$.
Recall that the parameter of the target distribution is $0.5$.

 \begin{figure}[ht]
  \includegraphics[trim={5cm 13cm 5cm 12.5cm},clip,width=0.80\textwidth]{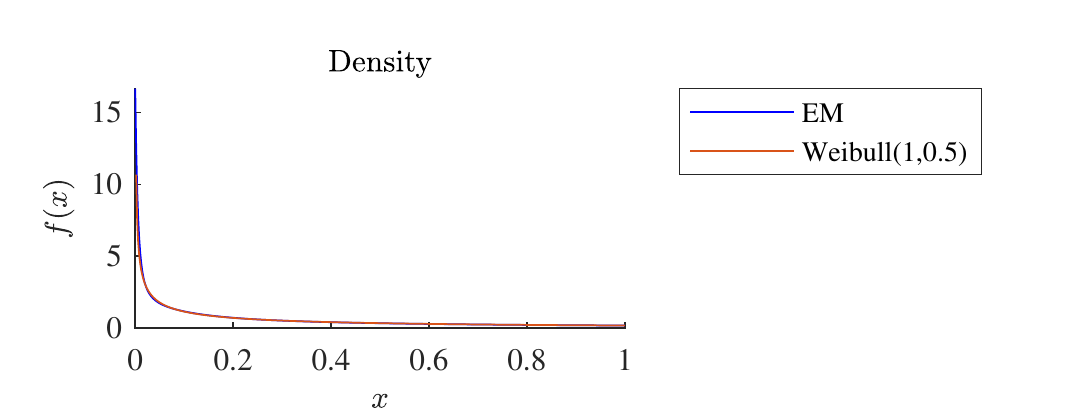}
  \caption{Density of the Weibull distribution with parameters $(1,0.5)$ and an approximation via 
   NPH models with EM adjustment.  
  }
  \label{Weibull.Dens}
 \end{figure} 

  \begin{figure}[ht]
  \includegraphics[trim={5cm 13cm 5cm 12.5cm},clip,width=0.70\textwidth]{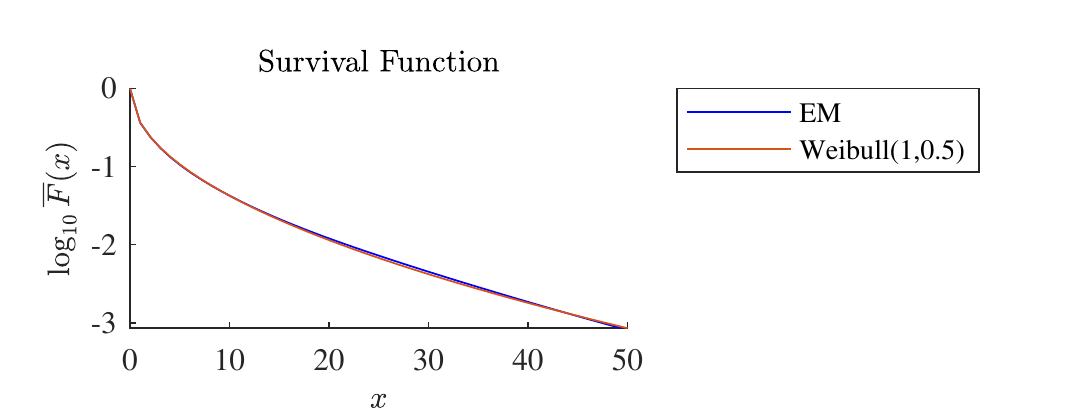}
  \includegraphics[trim={5cm 13cm 5cm 12.5cm},clip,width=0.70\textwidth]{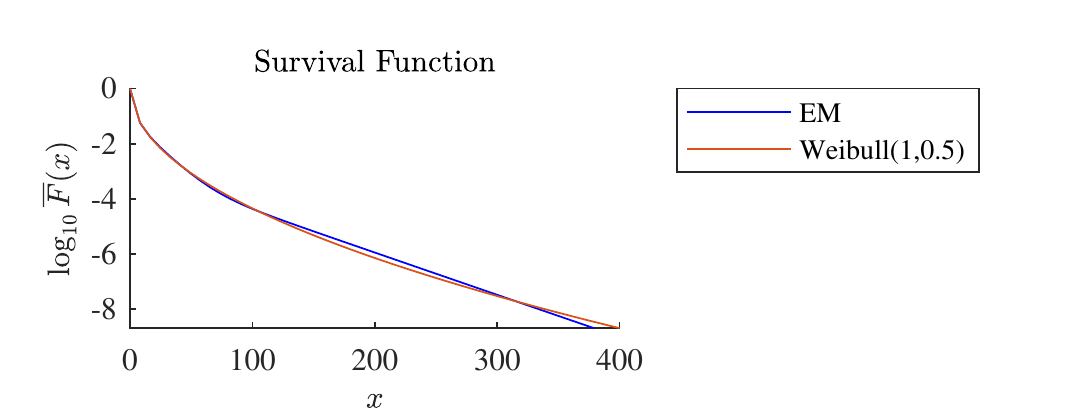}
  \caption{Survival functions of the $\mbox{Weibull}(1,0.5)$ and an approximation via an NDPH model with EM adjustment.
  The EM model is Weibullian with shape parameter $p= \widehat p(1+\widehat p)\approx 0.5386$.}  \label{Weibull.Survival1}
 \end{figure}

\end{Ex}

\begin{Ex}[Fitting to a theoretical distribution: Lognormal]\rm
Finally we consider a Lognormal distribution
with location parameter $\mu=0$ and dispersion parameter $\sigma^2=1$.
The lognormal distribution is heavy--tailed but its survival function ultimately
decays faster than a Regularly Varying function (but slower than a Weibullian function). 
We consider the lognormal case more difficult because a scaling distribution $\vect\pi(\vect\theta)$
for which the NPH model has the same tail behavior as the lognormal distribution is unknown.

Therefore, we consider two alternative  NPH models to adjust the data.  For both we have
selected a phase--type part has eight phases
and the scaling distribution is chosen as a discretization of a Lognormal distribution 
$\mbox{LN}(\mu,\sigma^2)$ and supported over the set $\{s_i=e^{i}:i=0,1,\dots\}$.  
The difference between the two models is that for the first one we let the EM
algorithm to estimate the values of the parameters $\mu$ and $\sigma$, while for
the second one we take these values to be fixed and equal to $0$ and $1$ respectively.
%
%
The results for the first model are given below
\begin{eqnarray*}
\hat{\mu}&=&0.1909, \ \ 
\hat{\sigma}=0.5979\\
\hat{\vect{\alpha}}&=&(  0.0000,    0.0000,    0.0000,       0,    0.0000,    0.0439,    0.1360,  0.8201)\\
\hat{\mat{T}}&=&
\begin{pmatrix}
   -7.1741 &   0.1557 &   0.3338 &   6.2239 &   0.0195 &   0.1985 &   0.1061 & 0.1366\\
    0.2056 &  -1.1622 &   0.0528 &   0.2434 &   0.1133 &   0.0519 &   0.0328 & 0.0782\\
    0.7528 &   0.0891 &  -6.2221 &   4.7969 &   0.1873 &   0.2546 &   0.1238 & 0.0176\\
    0.3840 &   0.1133 &   0.1303 &  -9.4035 &   0.1845 &   0.0654 &   0.2161 & 0.1597\\
    0.9550 &   0.4274 &   1.5430 &   0.9113 &  -5.2223 &   0.2926 &   1.0875 & 0.0056\\
    0.2341 &   0.0842 &   2.3467 &   0.0000 &   1.3606 &  -4.6381 &   0.0621 & 0.5504\\
    1.8359 &   0.0083 &   1.0178 &   0.0000 &   0.0211 &   0.2123 &  -3.3825 & 0.2870\\
    4.9757 &   0.0000 &   0.9689 &   0.0000 &   3.4702 &   0.2689 &   0.3086 &-9.9923
\end{pmatrix}
\end{eqnarray*}  

 \begin{figure}[ht]
  \includegraphics[trim={5cm 13cm 5cm 12.5cm},clip,width=0.80\textwidth]{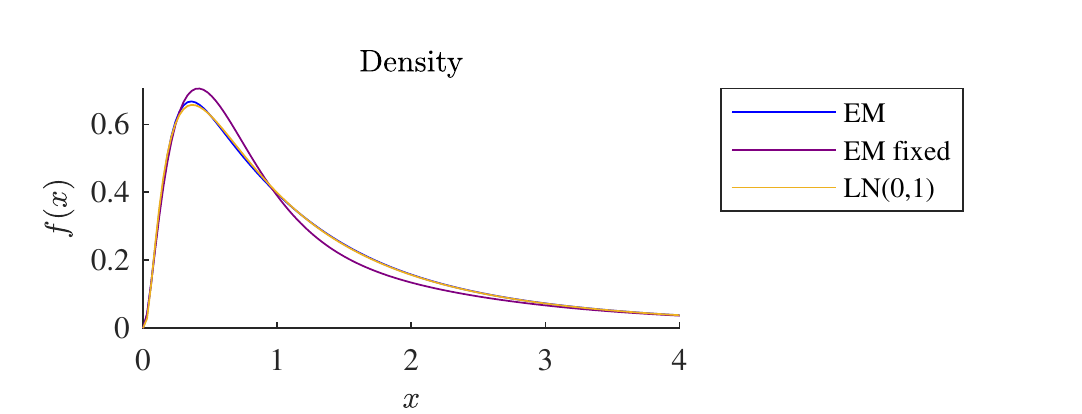}
  \caption{Density of the Lognormal distribution with parameters $(0,1)$ and an approximation via 
   NPH models with EM adjustment.}
  \label{Lognormal.Dens}
 \end{figure} 

  \begin{figure}[ht]
  \includegraphics[trim={5cm 13cm 5cm 12.5cm},clip,width=0.70\textwidth]{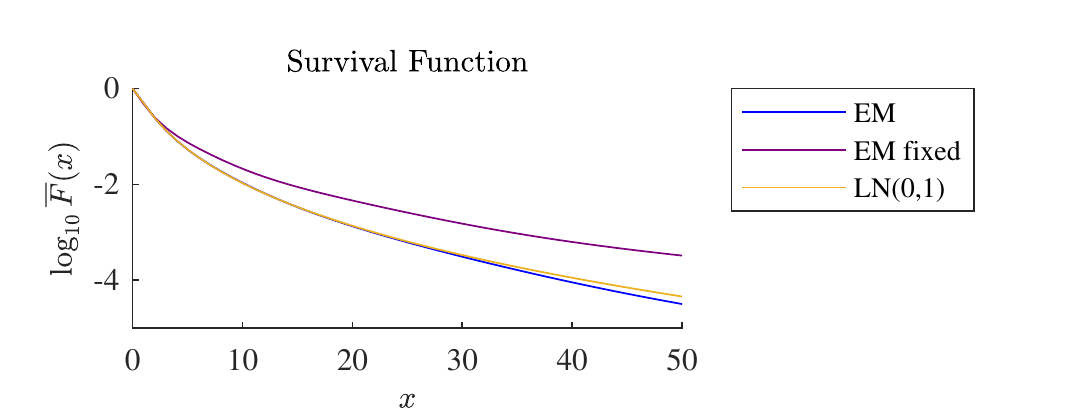}
  \includegraphics[trim={5cm 13cm 5cm 12.5cm},clip,width=0.70\textwidth]{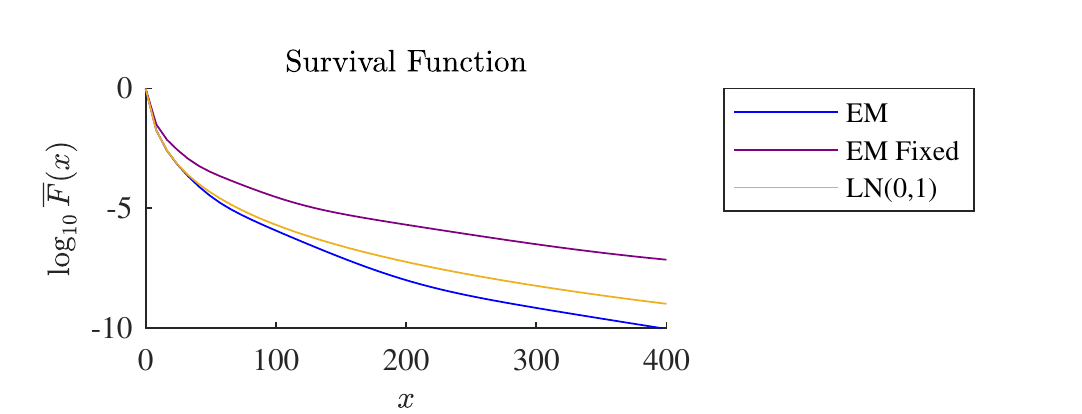}
  \caption{Survival functions of the Lognormal distribution 
   $\mbox{LN}(0,1)$ and an approximation via an NPH model with EM adjustment.}  \label{Lognormal.Survival1}
 \end{figure} 

From Figure \ref{Lognormal.Dens} we can observe that the adjustment of the body of the distributions of
the first model is excellent, but the tail probabilities differ.  In the lognormal case, the \emph{heaviness} of 
the distribution is mostly determined by the parameter $\sigma$.  The larger the parameter
$\sigma$ the more heavier its tail.  In our estimations we obtained an estimate
$\widehat\sigma=0.5979$ which suggest a lighter tail, and this is confirmed in the last panel of Figure \ref{Lognormal.Survival1}.

For our second model we obtained a very poor fit, but this is somewhat expected 
since the tail distribution  of the NPH model is very different from the
tail behavior of its scaling distribution (as
in the case of Regularly Varying or Weibullian cases).  In fact, \cite{RojasXie}
demonstrate that the tail probability of the NPH model is significantly  heavier
and different to the scaling distribution. More precisely, it is shown that
if $N\sim\mbox{LN}(0,1)$, then 
\begin{align*}
 \lim_{t\to\infty}\frac{\Prob(\tau N>t)}{\Prob(N>t)}=\infty.
\end{align*}  

\end{Ex}

\section{Conclusions}
{In this paper we have introduced a new methodology for adjusting phase--type scale mixture distributions
to heavy--tailed data.  While the model only requires the specification of the dimension of the underlying phase--type distribution and a parametric family of discrete scaling distributions, the suggested algorithm simultaneously fit  both the body and the tail of general distributions.}

{Since the class of NPH distributions are generally genuinely heavy tailed (if $N$ has unbounded support) and dense in the class of heavy tailed distributions with support on $\mathds{R}_+$, we may in principle approximate any heavy tailed distribution (data) arbitrarily close by a NPH distribution. In particular, for the case of regular varying and Weibullian distributions the aforementioned approximation is not only in the limit (denseness) but can be effectively carried out in praxis. }

\bibliographystyle{abbrv}

\end{document}